\documentclass[letterpaper, 10 pt, conference]{ieeeconf}
\IEEEoverridecommandlockouts
\IEEEoverridecommandlockouts
\usepackage{cite}
\usepackage{amsmath,amssymb,amsfonts}
\usepackage{booktabs}
\usepackage{hyperref}
\usepackage{cleveref}
\usepackage[table]{xcolor}

\usepackage{algorithm}
\usepackage{algpseudocode}

\usepackage{graphicx}
\usepackage{textcomp}
\newtheorem{theorem}{Theorem}
\newtheorem{lemma}{Lemma}
\newtheorem{remark}{Remark}
\newtheorem{definition}{Definition}
\newtheorem{assumption}{Assumption}

\newtheorem{proposition}{Proposition}
\newcommand{\wbar}{\bar}

\newcommand{\what}{\hat}
\usepackage{mathtools}

\title{ Distributed Online Convex Optimization with Nonseparable Costs and Constraints}

\author{
\authorblockN{Zhaoye Pan\textsuperscript{1}, Haozhe Lei\textsuperscript{2},~\IEEEmembership{Graduate Student Member, IEEE}, Fan Zuo\textsuperscript{3},~\IEEEmembership{Member, IEEE}, \\ Zilin Bian\textsuperscript{4},~\IEEEmembership{Member, IEEE}, and Tao Li\textsuperscript{1,$\dagger$},~\IEEEmembership{Member, IEEE}}
\thanks{\textsuperscript{1}Z. Pan and T. Li are with the Department of Systems Engineering, City University of Hong Kong, Hong Kong SAR, 999077, China (e-mail: iszhaoyepan@gmail.com, li.tao@cityu.edu.hk).
\textsuperscript{2}H. Lei is with the Department of Electrical and Computer Engineering, New York University, NY, 11201, USA (email: hl4155@nyu.edu). 
\textsuperscript{3}F. Zuo is with the Department of Civil and Urban Engineering, New York University, NY, 11201, USA (email: fz380@nyu.edu).
\textsuperscript{4}Z. Bian is with the Department of Civil Engineering Technology and Environmental Management Safety, Rochester Institute of Technology, NY, 14623, USA (email: zxbite@rit.edu).
\textsuperscript{$\dagger$}Corresponding author: T. Li (email: li.tao@cityu.edu.hk).
}
}

\begin{document}

\maketitle
\begin{abstract}
This paper studies distributed online convex optimization with time-varying coupled constraints, motivated by distributed online control in network systems. Most prior work assumes a separability condition: the global objective and coupled constraint functions are sums of local costs and individual constraints. In contrast, we study a group of agents, networked via a communication graph, that collectively select actions to minimize a sequence of \textit{nonseparable} global cost functions and to satisfy \textit{nonseparable} long-term constraints based on full-information feedback and intra-agent communication. We propose a distributed online primal-dual belief consensus algorithm, where each agent maintains and updates a local belief of the global collective decisions, which are repeatedly exchanged with neighboring agents. Unlike the previous consensus primal-dual algorithms under separability that ask agents to only communicate their local decisions, our belief-sharing protocol eliminates coupling between the primal consensus disagreement and the dual constraint violation, yielding sublinear regret and cumulative constraint violation (\textsc{ccv}) bounds, both in $O(\sqrt{T})$, where $T$ denotes the time horizon. Such a result breaks the long-standing $O(T^{3/4})$ barrier for \textsc{ccv} and matches the lower bound of online constrained convex optimization, indicating the online learning efficiency at the cost of communication overhead. 
\end{abstract}

\section{Introduction}
Distributed online convex optimization (\textsc{d-oco}) has emerged as a standard framework for multi-agent sequential decision-making under uncertainty, with applications spanning communication networks, power grids, and robotic systems \cite{lina23doco-review, tao22confluence}. A defining feature of \textsc{d-oco} is that the objective functions are revealed sequentially, and agents must take actions without prior knowledge of future costs, aiming to close the performance gap, known as the regret \cite{zinkevich03oco}, between the executed actions and the optimal in hindsight.

Starting from Zinkevich's seminal work on single-agent unconstrained \textsc{oco} \cite{zinkevich03oco}, later development has explored unconstrained \textsc{d-oco} \cite{zhang24doco}, constrained \textsc{oco} \cite{tianbao12const-constraint}, \textsc{oco} with time-varying constraints \cite{cao18oco-constraint, muthirayan23oco-pred, sinha24opt-coco}, and most recently, \textsc{d-oco} with time-varying constraints, where agents' joint actions are subject to a constraint function changing over time that is previously unknown and sequentially revealed to agents \cite{yi19doco, yangtao20doco-bandit, yuan20ma-coco}. Such a constrained formulation is motivated by the operational and safety constraints in the multi-agent system. 

Most existing \textsc{d-oco} formulations assume that both the objective and constraint functions are \emph{separable}, i.e., expressed as sums of local agent-specific functions \cite{yuan20ma-coco}. Under such a separability condition, the seminal online gradient descent (\textsc{ogd}) \cite{zinkevich03oco} leads to a consensus algorithm where agent share with their neighbors the local gradient in both primal and dual variables with respect to the Lagrangian \cite{zavlanos16doco}. 

However, the separability assumption fails to capture a broad class of online distributed control problems featuring team-based decision-making \cite{yuksel24control-is}, where system performance and safety metrics depend on the \emph{joint configuration} of all agents. One prominent example is the recent development in risk-aware cooperative traffic sensing \cite{tao24picol,tao25dima}, where surveillance cameras are sequentially tilted in different directions to collectively maximize the captured urban-level traffic volume, which cannot be decomposed into camera-specific objectives. Essentially, such nonseparability points to the anonymity (agent-agnostic) and symmetry characteristics commonly observed in large-population multi-agent decision-making, with fruitful applications in traffic assignment \cite{pan-tao22noneq, pan-tao23delay} and swarm robotics \cite{lauri23pomdp-robot}, and broader relevance to decentralized Markov decision processes \cite{lauri23pomdp-robot} and generalized Nash equilibrium in potential games \cite{lina13state_game,guanze21local-gne, shutian23erm}.

Extending the state-of-the-art (\textsc{sota}) consensus algorithms \cite{yi19doco, yangtao20doco-bandit, yuan20ma-coco} to the nonseparable scenarios is not straightforward. Due to separability, the global objective's partial derivative with respect to each agent's action equals the gradient of that agent's local objective. By sharing their actions, agents effectively communicate their gradient updates, whose aggregation gives the global gradient update. In contrast, nonseparability makes it difficult for each agent to compute the gradient without knowing the joint actions.  

Our intuition for tackling this challenge is simple yet effective: if agents need to know the joint actions for computing gradients, they can maintain local beliefs about joint actions, which are then calibrated through communication and local gradient update. Towards this end, we develop a \emph{distributed online primal--dual algorithm with belief consensus} (\textsc{dopbc}), where each agent maintains a belief matrix with each column representing its estimate of others' actions to be exchanged among neighbors. \textsc{dopbc} shares with early works on distributed optimization with nonseparable global objectives \cite{lina13state_game, raginsky16nonseparable} a similar spirit in exchanging agents' beliefs of joint actions. Yet, \textsc{dopbc} targets the constrained case in which its belief consensus includes both primal and dual estimates under weighted averaging, yielding a simpler communication protocol than the prior works under which we find that the constraint violation can be reduced.

\paragraph{Contribution}
Our contribution is two-fold. We first develop a belief-sharing scheme in a distributed online primal-dual algorithm to address the nonseparability challenge in \textsc{d-oco}. Additionally, we establish tight upper bounds on the regret and \textsc{ccv} growth of the proposed \textsc{dopbc} with refined dependence on the number of agents and communication network connectivity. Specifically, we show in \Cref{thm:violation} and \ref{thm:static_regret} that both regret and \textsc{ccv} are of $O(\sqrt{T})$ in a $T$-horizon \textsc{d-oco}, breaking the existing $O(T^{3/4})$ barrier for \textsc{ccv} \cite{yi19doco, yangtao20doco-bandit, yuan20ma-coco} and matching the lower bound of $O(\sqrt{T})$ for both regret and \textsc{ccv} in \textsc{oco} with time-varying constraints \cite{sinha24opt-coco} at large. In other words, one cannot further improve the regret or \textsc{ccv} bound without sacrificing the other. \Cref{cor:global_belief_advantage} further unveils the advantage of sharing beliefs over conventional consensus methods that share local information: belief-sharing eliminates the coupling between the
primal consensus disagreement and the dual constraint violation, justifying \textsc{dopbc}'s improved bounds.  

\paragraph{Related Works} Most existing online learning algorithms in \textsc{oco} are based on the seminal \textsc{ogd} and its mirror descent variants \cite{zinkevich03oco, zhang24doco, pan24mirror}. Moving to \textsc{d-oco}, \textsc{ogd} leads to distributed primal-dual methods resulting from Lagrangians, which yield a trade-off between regret growth and \textsc{ccv} growth. \Cref{tab:summary} summarizes the \textsc{sota} results \textsc{d-oco} with time-varying convex costs and constraints. As shown in \Cref{tab:summary}, when the trade-off parameter $c=1/2$, one recovers the optimal regret bound in \cite{yi19doco, yangtao20doco-bandit, yuan20ma-coco}, albeit with greater \textsc{ccv} of $O(T^{3/4})$. Our belief consensus achieves a tight \textsc{ccv} (when $c=1/2$) bound without affecting the regret at the cost of sharing an entire belief matrix rather than decision vectors as in conventional consensus algorithms.
\begin{table}[!ht]
\vspace{-1em}
    \centering
        \caption{A summary of existing regret and \textsc{ccv} bounds. $c\in (0,1)$ denotes the trade-off parameter. }
        \resizebox{\linewidth}{!}{
    \begin{tabular}{llll}
    \toprule
      Refs   & Regret & \textsc{ccv} & Comment \\
      \midrule
      \cite{yi19doco}   &  $O(T^{\max\{1-c,c\}})$ & $O(T^{1-c/2})$ & consensus, full information\\
      \cite{yangtao20doco-bandit} & $O(T^{\max\{1-c,c\}})$ & $O(T^{1-c/2})$ & consensus, bandit feedback\\
      \cite{yuan20ma-coco} & $O(T^{\max\{1-c,c\}})$ & $O(T^{1-c/2})$ & consensus, full information\\
      \cite{sinha24opt-coco} & $O(\sqrt{T})$ & $O(\sqrt{T}\ln T)$ & one agent, full information\\
      Ours & $O({T}^{\max\{1-c,c\}})$ & $O({T}^{\max\{1-c,c\}})$ & belief, full information\\
      \rowcolor{gray!20} \cite{sinha24opt-coco}  & $O(\sqrt{T})$ & $O(\sqrt{T})$ & \textsc{oco} lower bounds\\
      \bottomrule
    \end{tabular}
    }
    \vspace{-1em}
    \label{tab:summary}
\end{table}

\section{Problem Formulation}

Consider $N$ agents indexed by the set $\mathcal{V}=\{1,\ldots,N\}$ networked by an undirected graph $G=(\mathcal{V},\mathcal{E})$. At each time $t\in\{1,\ldots,T\}$, agent $i\in\mathcal{V}$ selects a local action $\mathbf{a}_{i,t}\in\mathcal{X}_i\subset\mathbb{R}^{d_i}$; $\mathcal{X}_i$ is convex and compact. The joint action is $\mathbf{a}_t \triangleq (\mathbf{a}_{1,t}^\top,\ldots,\mathbf{a}_{N,t}^\top)^\top \in \mathcal{X}\triangleq\prod_{i=1}^N \mathcal{X}_i \subset \mathbb{R}^{d}$, $d\triangleq\sum_{i=1}^N d_i$.

After executing $\mathbf{a}_t$, the environment reveals a global and nonseparable cost function $f_t:\mathcal{X}\to\mathbb{R}$ and a vector of $m$ nonseparable constraint functions $\mathbf{g}_t:\mathcal{X}\to\mathbb{R}^m$, whose $k$-th component is denoted by $g_{k,t}:\mathcal{X}\to\mathbb{R}$. The  objective is to generate a sequence of joint actions $\{\mathbf{a}_t\}_{t=1}^T$ that minimizes the cumulative cost while satisfying the constraints in the long run: denoting $\preceq$ element-wise inequality,
\begin{equation}\label{eq:problem_formulation}
\begin{aligned}
\min_{\{\mathbf{a}_t\}_{t=1}^T \subseteq \mathcal{X}} \quad & \sum_{t=1}^T f_t(\mathbf{a}_t) \\
\text{s.t.}\quad & \mathbf{g}_t(\mathbf{a}_t)\preceq \mathbf{0}, \quad \forall t\in\{1,\ldots,T\},
\end{aligned}
\end{equation}

We evaluate the performance of a distributed online algorithm via static regret and cumulative constraint violation.

\begin{definition}[Static Regret]
Let $\mathbf{a}^\star$ be the best fixed feasible joint action in hindsight:
\begin{equation*}
\mathbf{a}^\star \in \arg\min_{\mathbf{a}\in\mathcal{X}} \left\{ \sum_{t=1}^T f_t(\mathbf{a}) \ \middle|\  \mathbf{g}_t(\mathbf{a})\preceq \mathbf{0},\ \forall t\in\{1,\ldots,T\} \right\}.
\end{equation*}
The regret is  $\mathcal{R}_S(T) \triangleq \sum_{t=1}^T f_t(\mathbf{a}_t) - \sum_{t=1}^T f_t(\mathbf{a}^\star)$.
\end{definition}

\begin{definition}[Cumulative Constraint Violation]
For $k\in\{1,\ldots,m\}$, the cumulative constraint violation is
$\mathcal{C}_k(T) \triangleq \sum_{t=1}^T \big[g_{k,t}(\mathbf{a}_t)\big]^+$,
where $[z]^+\triangleq\max\{0,z\}$. 
\end{definition}

\section{Distributed Online Primal-Dual with Belief Consensus}

To solve the cooperative online optimization problem in \eqref{eq:problem_formulation} in the sense that $\mathcal{R}_S(T)$ and $\mathcal{C}_k(T)$ grow sublinearly with $T$ for all $k$, we propose a distributed algorithm: \emph{Distributed Online Primal-Dual with Belief Consensus} (\textsc{dopbc}). The algorithm is built upon a primal--dual framework applied to the time-varying Lagrangian
\begin{equation*}
\mathcal{L}_t(\mathbf{x}, \mathbf{\lambda}) \triangleq f_t(\mathbf{x}) + \mathbf{\lambda}^\top \mathbf{g}_t(\mathbf{x}),
\end{equation*}
where $\mathbf{x}\in\mathcal{X}$ denotes the global joint action and $\mathbf{\lambda}\in\mathbb{R}^m_{\ge 0}$ is the Lagrange multiplier vector with respect to constraints.

Each agent $i$ maintains a local estimate of the global primal variable $\mathbf{x}_{i,t}\in\mathcal{X}$ (which shall be distinguished from the local action $\mathbf{a}_{i,t}\in \mathcal{X}_i$) and the global dual variable $\mathbf{\lambda}_{i,t}\in\mathbb{R}^m_{\ge 0}$. The algorithm follows a \emph{consensus-then-update} structure: agents first fuse local estimates through neighbor communication, then form and execute local actions, and finally update their primal and dual estimates based on the observed global feedback.

A key challenge is that agent $i$ controls only its own action block and therefore cannot compute the full gradient of the Lagrangian. To address this, each agent computes a \emph{primal pseudo-gradient}, defined as the partial gradient of the Lagrangian with respect to its own decision, evaluated at the executed joint action and estimated dual variables (see \eqref{eq:pseudo_grad}) while treating other agents' actions as fixed. This pseudo-gradient captures the marginal effect of agent $i$'s action on the global objective and constraints and enables fully decentralized primal updates. The complete procedure is summarized in Algorithm~\ref{alg:dopbc}. 

\begin{algorithm}[t]
\caption{\textsc{dopbc}: Distributed Online Primal-Dual with Belief Consensus}
\label{alg:dopbc}
\begin{algorithmic}[1]
\State \textbf{Input:} Step-size sequence $\{\alpha_t\}_{t>0}$; a doubly stochastic matrix $W$ compatible with graph $G$ (Assumption~\ref{ass:net}).
\State \textbf{Initialization:} Each agent $i\in\mathcal{V}$ initializes the belief $\mathbf{x}_{i,1}\in\mathcal{X}$ and $\mathbf{\lambda}_{i,1}=\mathbf{0}\in\mathbb{R}^m_{\ge 0}$.

\For{$t = 1,2,\ldots,T$}
    \ForAll{agents $i\in\mathcal{V}$ \textbf{in parallel}}
    \Comment{\texttt{Belief consensus}} 
    \State Agent $i$ exchanges $(\mathbf{x}_{i,t},\mathbf{\lambda}_{i,t})$ with neighbors $j\in\mathcal{N}_i$ and computes $$\widehat{\mathbf{x}}_{i,t} = \sum_{j=1}^N W_{ij}\mathbf{x}_{j,t},\quad \widehat{\mathbf{\lambda}}_{i,t} = \sum_{j=1}^N W_{ij}\mathbf{\lambda}_{j,t}.$$

    \State Agent $i$ forms its local action $\mathbf{a}_{i,t} =[\widehat{\mathbf{x}}_{i,t}]_i$ (under projection)
    and the joint action $\mathbf{a}_t=(\mathbf{a}_{1,t}^\top,\ldots,\mathbf{a}_{N,t}^\top)^\top$ is executed. All agents observe $f_t(\cdot)$ and $\mathbf{g}_t(\cdot)$.
    \Comment{\texttt{Primal update}}
    \State  Agent $i$ computes the primal pseudo-gradient
    \begin{equation}\label{eq:pseudo_grad}
        \mathbf{d}_{i,t} = \nabla_i f_t(\widehat{\mathbf{x}}_{i,t}) + \sum_{k=1}^m [\widehat{\mathbf{\lambda}}_{i,t}]_k \nabla_i g_{k,t}(\widehat{\mathbf{x}}_{i,t}),
    \end{equation}
    and updates its global estimate block-wise as
    \begin{equation}\label{eq:primal_update_blockwise}
        [\mathbf{x}_{i,t+1}]_j =
        \begin{cases}
        \mathcal{P}_{\mathcal{X}_i}\!\left([\widehat{\mathbf{x}}_{i,t}]_i - \alpha_t \mathbf{d}_{i,t}\right), & j=i, \\
        [\widehat{\mathbf{x}}_{i,t}]_j, & j\neq i.
        \end{cases}
    \end{equation}
    \Comment{\texttt{Dual update}}
    \State Agent $i$ updates the dual estimate
    \begin{equation}\label{eq:dual_update}
        \mathbf{\lambda}_{i,t+1}
        = \mathcal{P}_{\mathbb{R}^m_{\ge 0}}\!\left(\widehat{\mathbf{\lambda}}_{i,t} + \alpha_t \mathbf{g}_t(\widehat{\mathbf{x}}_{i,t})\right).
    \end{equation}
    
    \EndFor
\EndFor
\end{algorithmic}
\end{algorithm}

\paragraph{Assumptions}
Our analysis rests on standard assumptions on the sets and functions regularity, and network topology, commonly seen in prior works \cite{yi19doco, yangtao20doco-bandit, yuan20ma-coco}.
\begin{assumption}[Regularity]\label{ass:std}
For each time step $t$, the following hold:
\begin{enumerate}
    \item \textbf{Convexity and Compactness:} The cost $f_t:\mathcal{X}\to\mathbb{R}$ and each constraint component $g_{k,t}:\mathcal{X}\to\mathbb{R}$, $k\in\{1,\ldots,m\}$, are convex. The feasible set $\mathcal{X}\subseteq\mathbb{R}^d$ is convex and compact with diameter bounded by $D$, i.e., $\|\mathbf{x}-\mathbf{y}\|\le D$ for all $\mathbf{x},\mathbf{y}\in\mathcal{X}$.
    \item \textbf{Bounded Gradients:} The functions have uniformly bounded gradients on $\mathcal{X}$. There exist constants $G_f>0$ and $G_g>0$ such that for all $\mathbf{x}\in\mathcal{X}$ and all $t$, $\|\nabla f_t(\mathbf{x})\|\le G_f$ and $\|\nabla \mathbf{g}_t(\mathbf{x})\|_F \le G_g$,
    where $\nabla \mathbf{g}_t(\mathbf{x})\in\mathbb{R}^{m\times d}$ denotes the Jacobian of $\mathbf{g}_t$ and $\|\cdot\|_F$ is the Frobenius norm.
    \item \textbf{Bounded Dual Domain:} The dual variable is projected onto the compact set $\Lambda \triangleq [0,\Lambda_{\max}]^m$. Consequently, $\|\mathbf{\lambda}\|\le \Lambda_{\mathrm{bound}}\triangleq\sqrt{m}\Lambda_{\max}$ for all $\mathbf{\lambda}\in\Lambda$.
\end{enumerate}
\end{assumption}

\begin{assumption}[Network Connectivity]\label{ass:net}
The communication graph $G=(\mathcal{V},\mathcal{E})$ is fixed, undirected, and connected. The consensus matrix $W=[W_{ij}]\in\mathbb{R}_{\geq 0}^{N\times N}$ is compatible with $G$ in the sense that $W_{ij}=0$ for $i\neq j$ and $(i,j)\notin\mathcal{E}$, and it satisfies:
\begin{enumerate}
    \item \textbf{Doubly stochasticity:}  $W\mathbf{1}=\mathbf{1}$ and $\mathbf{1}^\top W=\mathbf{1}^\top$.
    \item \textbf{Spectral property:} 
    $W$ is positive semi-definite, and its second largest eigenvalue, $\sigma$, satisfies $\sigma\in [0, 1)$.
\end{enumerate}
\end{assumption}

\begin{remark}[Justification for Bounded Dual Domain]
\label{rem:dual_bound_justification}
The assumption of a bounded dual domain is well-founded in the literature of distributed constrained optimization. For instance, works employing Augmented Lagrangian methods \cite{yuan20ma-coco} have remarked that the penalty term inherently serves to upper bound the dual variables. More fundamentally, Slater's condition guarantees that the optimal Lagrange multipliers satisfy a finite norm bound \cite{yi19doco}. Consequently, projecting dual variables onto a sufficiently large compact set is a standard practice adopted in prior works (e.g., \cite{zavlanos16doco}) to ensure algorithm stability without loss of generality.
\end{remark}

\subsection{Theoretic Analysis}
This subsection first presents several key lemmas that form the basis of the regret and constraint violation analysis. We begin by defining the network-wide averages of the primal and dual estimates, along with the associated total primal consensus error:
\begin{align*}
\wbar{\mathbf{x}}_t \triangleq \frac{1}{N}\sum_{i=1}^N \mathbf{x}_{i,t}, 
\wbar{\boldsymbol{\lambda}}_t \triangleq \frac{1}{N}\sum_{i=1}^N \boldsymbol{\lambda}_{i,t}, 
\delta_{x,t}^2 \triangleq \sum_{i=1}^N \|\mathbf{x}_{i,t} - \wbar{\mathbf{x}}_t\|^2.
\end{align*}

The theoretical analysis effectively handles the distributed coupled constraints through two fundamental steps. Lemma \ref{lemma:consensus_recur} quantifies the primal disagreement among agents by bounding the consensus error, ensuring that local estimates track the network average. Lemma \ref{lemma:drift} establishes the optimization progress via a one-step Lyapunov drift analysis. This decomposes the instantaneous Lagrangian difference into a telescoping term and an error term dependent on the consensus discrepancy. These two lemmas provide the necessary building blocks for the global convergence analysis.

\begin{lemma}[One-step Primal Consensus Error]\label{lemma:consensus_recur}
Under Assumptions~\ref{ass:std} and~\ref{ass:net}, for any constant $\beta>0$, the primal consensus error $\delta_{x,t}^2$ generated by Algorithm~\ref{alg:dopbc} satisfies
$
    \delta_{x,t+1}^2
    \le (1+\beta)\sigma^2 \delta_{x,t}^2
    + (1+1/\beta)\alpha^2 N G_d^2,
$
where $G_d \triangleq G_f + \Lambda_{\max} G_g$.
\end{lemma}
\begin{lemma}[One-Step Primal-Dual Drift Analysis]\label{lemma:drift}
Under Assumptions~\ref{ass:std} and~\ref{ass:net} and the constant step size $\alpha_t=T^{-c}$, $c\in (0,1)$
for any comparator $(\boldsymbol{x}^*,\boldsymbol{\lambda}^*)$ with $\boldsymbol{\lambda}^*\in[0,\Lambda_{\max}]^m$, the following inequality holds for each $t=1,\ldots,T$:
\begin{align}
\label{eq:lemma_drift_fixed}
&f_t(\widehat{\mathbf{x}}_t) - f_t(\mathbf{x}^*)
 + \bigl\langle \wbar{\boldsymbol{\lambda}}_t - \boldsymbol{\lambda}^*,\,
   \mathbf{g}_t(\widehat{\mathbf{x}}_t) \bigr\rangle \nonumber \\ 
&\le
  \frac{1}{2\alpha}
  \Bigl(
  \|\wbar{\mathbf{x}}_t - \mathbf{x}^*\|^2
  -\|\wbar{\mathbf{x}}_{t+1} - \mathbf{x}^*\|^2
  \Bigr) \nonumber \\
&
  +\frac{1}{2\alpha}
  \Bigl(
  \|\wbar{\boldsymbol{\lambda}}_t - \boldsymbol{\lambda}^*\|^2
  -\|\wbar{\boldsymbol{\lambda}}_{t+1} - \boldsymbol{\lambda}^*\|^2
  \Bigr) +\frac{\alpha}{2}K_c + C_E\delta_{x,t}.\nonumber
\end{align}
where $K_c$ is a constant bounding the norm of the primal pseudo-gradients, and $C_E$ is a positive constant depending only on the Lipschitz parameters $(G_f,G_g,\Lambda_{\max})$,
\end{lemma}

The main results are derived by aggregating the one-step bounds over the entire time horizon $T$. Specifically, the proof of Theorem \ref{thm:violation} utilizes the dual drift component from Lemma \ref{lemma:drift} to bound the accumulated constraint violation. Subsequently, for Theorem \ref{thm:static_regret}, we substitute this violation bound and the consensus error bound (from Lemma \ref{lemma:consensus_recur}) back into the primal drift decomposition. 
\begin{theorem}[Cumulative Constraint Violation]\label{thm:violation}
Given Assumptions~\ref{ass:std} and~\ref{ass:net}, let the step size be set as $\alpha_t = T^{-c}$ where $c \in (0,1)$.
For the sequence of global average primal estimates $\{\widehat{\mathbf{x}}_t\}_{t=1}^T$ generated by Algorithm~\ref{alg:dopbc}, the \textsc{ccv} of each component $k\in\{1,\ldots,m\}$ satisfies
\begin{align*}
\mathcal{C}_k(T)
&\coloneqq
\sum_{t=1}^T \big[g_t^{(k)}(\widehat{\mathbf{x}}_t)\big]^+ \nonumber \\
\;&\le\;
\Lambda_{\max} T^c
+
\left( \frac{L_g\,\sqrt{N}\,G_d}{1-\sigma} \right) T^{1-c} \nonumber \\
& = O\!\left( T^c + \frac{\sqrt{N}}{1-\sigma}T^{1-c} \right),\quad G_d \triangleq G_f + \Lambda_{\max}G_g.
\end{align*}
\end{theorem}

\begin{theorem}[Static Regret]\label{thm:static_regret}
Given Assumptions~\ref{ass:std} and~\ref{ass:net}, let the step size be set as $\alpha_t = T^{-c}$ with constant $c \in (0,1)$.
For any fixed comparator $\mathbf{x}^* \in \mathcal{X}$, the static regret of Algorithm~\ref{alg:dopbc} satisfies
\begin{align*}
\mathcal{R}_S(T)
\;&\le\;
\frac{\hat{D}^2}{2} T^c 
+
\left( \frac{K_c+G_g^2}{2} + \frac{C_{\mathrm{net}}\sqrt{N}}{1-\sigma} \right) T^{1-c} \nonumber\\
& = O\!\left( T^c + \frac{\sqrt{N}}{1-\sigma}T^{1-c} \right),
\end{align*}
where $\hat{D}^2 \triangleq D^2+\Lambda_{\max}^2$, $K_c$ bounds the pseudo-gradient norm, and $C_{\mathrm{net}}$ depends on gradient bounds $(G_f,G_g,\Lambda_{\max})$.
Specifically, setting $c=1/2$ recovers the optimal $O(\sqrt{T})$.
\end{theorem}

\begin{proposition}[Primal-Dual Decoupling]
\label{cor:global_belief_advantage}
Consider a distributed online optimization problem with a separable objective and a separable coupled constraint
\[
f_t(x)=\sum_{i=1}^N f_{i,t}(x_i),\quad g_t(x)=\sum_{i=1}^N g_{i,t}(x_i)\le 0,
\]
under similar regularity assumptions as in our case. 

Under the conventional consensus framework, e.g., \cite{yuan20ma-coco}, where agents share local decisions $x_{i,t}\in \mathcal{X}_i$, the cumulative primal consensus error admits an upper bound of the form 
\begin{equation}\label{eq:coupling_bound}
\sum_{t=1}^T \sum_{i=1}^N \|x_{i,t}-\bar x_t\|
\;\le\;
C_1 \sum_{t=1}^T \alpha_t
+
C_2 \sum_{t=1}^T \alpha_t
\sum_{i=1}^N [g_{i,t}(x_{i,t})]_+ ,
\end{equation}
which reveals an explicit coupling between primal disagreement and constraint violations.
As a result, choosing step size $\alpha_t=T^{-1/2}$ ($c=1/2$) to ensure the optimal regret
$\mathrm{Regret}(T)=O(\sqrt{T})$ inevitably leads to $\sum_{t=1}^T \sum_{i=1}^N [g_{i,t}(x_{i,t})]_+
= O(T^{3/4})$, as shown in \cite{yuan20ma-coco}.

In contrast, under the proposed belief consensus, the primal consensus error among agents obeys a decoupled recursion
$
\delta_{x,t+1}
\le
\sigma\,\delta_{x,t} + \alpha_t G$,
$\delta_{x,t} \triangleq \frac{1}{N}\sum_{i=1}^N \|x_{i,t}-\bar x_t\|,
$
which is independent of constraint violations.
Consequently,
$
\sum_{t=1}^T \delta_{x,t}
=
O(\sum_{t=1}^T \alpha_t).
$
With the constant step size rule, the algorithm simultaneously achieves
\[
\mathrm{Regret}(T)=O(\sqrt{T}),
\qquad
\sum_{t=1}^T \sum_{i=1}^N [g_{i,t}(x_{i,t})]_+
=
O(\sqrt{T}).
\]
\end{proposition}
\begin{remark}
The gap between the $O(T^{3/4})$ violation under decision sharing and the
$O(\sqrt{T})$ violation achieved here is not due to step-size tuning,
but rather stems from a fundamental structural difference:
global belief sharing removes the intrinsic coupling between primal disagreement
and constraint violation, enabling optimal-order bounds for both metrics simultaneously.
\end{remark}
\section{Conclusion}
This paper investigated distributed online convex optimization with long-term coupled constraints under a global belief sharing framework.
We proposed a distributed online primal--dual belief consensus algorithm that enables agents
to cooperatively handle global objectives and constraints.
By leveraging global belief sharing, the proposed algorithm achieves $O(\sqrt{T})$ static regret
and $O(\sqrt{T})$ cumulative constraint violation with a fixed step size.
These guarantees demonstrate that effective learning and long-term feasibility can be attained
simultaneously in distributed online settings with globally coupled objectives and constraints.

Several extensions are worth future investigation.
First, it is of interest to extend the current analysis to time-varying or adaptive step-size schemes, which may further improve empirical performance and robustness.
Second, relaxing the full-information assumption to bandit or zeroth-order feedback remains an important direction, where the interaction between exploration, feasibility, and network disagreement poses new challenges.
\bibliographystyle{ieeetr}
\bibliography{ref}

\appendix

\subsection{Proof of Lemma \ref{lemma:consensus_recur}}
Define the intermediate unprojected variable $\mathbf{u}_{i,t+1} = \what{\mathbf{x}}_{i,t} - \alpha \mathbf{p}_{i,t}$, where $\mathbf{p}_{i,t}$ is the zero-padded local pseudo-gradient, whose $i$-th block is $\mathbf{d}_{i,t}$ and all other blocks are zero. Using the identity $\sum_{i=1}^N \|\mathbf{v}_i - \bar{\mathbf{v}}\|^2 = \frac{1}{2N}\sum_{i,j}\|\mathbf{v}_i - \mathbf{v}_j\|^2$ and the non-expansive property of projection, we have
\begin{align}
\delta_{x,t+1}^2 &= \sum_{i=1}^N \|\mathbf{x}_{i,t+1} - \wbar{\mathbf{x}}_{t+1}\|^2 \leq \sum_{i=1}^N \|\mathbf{u}_{i,t+1} - \wbar{\mathbf{u}}_{t+1}\|^2. \nonumber
\end{align}
Expanding $\mathbf{u}_{i,t+1} - \wbar{\mathbf{u}}_{t+1} = (\what{\mathbf{x}}_{i,t} - \wbar{\what{\mathbf{x}}}_t) - \alpha(\mathbf{p}_{i,t} - \wbar{\mathbf{p}}_t)$ and applying Young's inequality with parameter $\beta > 0$:
\begin{align}
\|\mathbf{u}_{i,t+1} - \wbar{\mathbf{u}}_{t+1}\|^2 &\leq (1+\beta)\|\what{\mathbf{x}}_{i,t} - \wbar{\what{\mathbf{x}}}_t\|^2 \nonumber \\
&\quad + (1+\beta^{-1})\alpha^2\|\mathbf{p}_{i,t} - \wbar{\mathbf{p}}_t\|^2. \nonumber
\end{align}
For the consensus term, using the spectral gap property of $W$:
$
\sum_{i=1}^N \|\what{\mathbf{x}}_{i,t} - \wbar{\what{\mathbf{x}}}_t\|^2 \leq \sigma^2 \delta_{x,t}^2,
$
where $\sigma = \|W - \frac{1}{N}\mathbf{1}\mathbf{1}^\top\| < 1$. 
For the gradient term, since $$\|\mathbf{p}_{i,t}\| \leq G_d, \quad 
\sum_{i=1}^N \|\mathbf{p}_{i,t} - \wbar{\mathbf{p}}_t\|^2 \leq N G_d^2.$$
Combining these bounds and choosing $\beta = \frac{1-\sigma}{1+\sigma}$ yields
\begin{equation}
\delta_{x,t+1}^2 \leq \sigma^2 (1+\beta) \delta_{x,t}^2 + (1+\beta^{-1}) N G_d^2 \alpha^2. \nonumber
\end{equation}
Simplifying with $\tilde{\sigma} = \sigma\sqrt{1+\beta} < 1$ completes the proof.
\QED

\subsection{Proof of Lemma \ref{lemma:drift}}
Define the instantaneous Lagrangian difference:
\[
\Delta_t \triangleq f_t(\widehat{\mathbf{x}}_t) - f_t(\mathbf{x}^*) + \langle \wbar{\boldsymbol{\lambda}}_t - \boldsymbol{\lambda}^*, \mathbf{g}_t(\widehat{\mathbf{x}}_t) \rangle.
\]

\paragraph*{Step 1 (Decomposition)} 
Define the average constraint vector effectively tracked by the dual update as $\wbar{\mathbf{g}}_t \coloneqq \frac{1}{N}\sum_{i=1}^N \mathbf{g}_{i,t}(\mathbf{x}_{i,t})$. By the convexity of $f_t$ and $\mathbf{g}_t$, we decompose $\Delta_t$ into an interpolated averaged drift term and a decentralization error $\mathcal{E}_t$:
\begin{equation} \label{eq:decomp}
\Delta_t \leq \langle \wbar{\mathbf{p}}_t, \wbar{\mathbf{x}}_t - \mathbf{x}^* \rangle + \langle \wbar{\mathbf{g}}_t, \wbar{\boldsymbol{\lambda}}_t - \boldsymbol{\lambda}^* \rangle + \mathcal{E}_t.
\end{equation}
The error term $\mathcal{E}_t$ captures the discrepancy between local variables and their averages, explicitly defined as:
\[
\mathcal{E}_t \triangleq \langle \nabla f_t(\wbar{\mathbf{x}}_t) - \wbar{\mathbf{p}}_t, \wbar{\mathbf{x}}_t - \mathbf{x}^* \rangle + \langle \wbar{\boldsymbol{\lambda}}_t - \boldsymbol{\lambda}^*, \mathbf{g}_t(\wbar{\mathbf{x}}_t) - \wbar{\mathbf{g}}_t \rangle.
\]
By Lipschitz continuity and the boundedness of the dual variable ($\|\wbar{\boldsymbol{\lambda}}_t\| \leq \Lambda_{\max}$), applying the Cauchy-Schwarz inequality yields the bound $|\mathcal{E}_t| \leq C_E \delta_{x,t}$, where $C_E$ is a constant depending on $(G_f, G_g, \Lambda_{\max})$.

\paragraph*{Step 2 (Primal drift)} By Jensen's inequality and non-expansiveness of projection:
\begin{align}
\langle \wbar{\mathbf{p}}_t, \wbar{\mathbf{x}}_t - \mathbf{x}^* \rangle &\leq \frac{1}{2\alpha} \left( \|\wbar{\mathbf{x}}_t - \mathbf{x}^*\|^2 - \|\wbar{\mathbf{x}}_{t+1} - \mathbf{x}^*\|^2 \right) \nonumber \\
&\quad + \frac{\alpha}{2} K_c + C\delta_{x,t}. \nonumber
\end{align}

\paragraph*{Step 3 (Dual drift)} Similarly, for the dual variable:
\begin{align}
\langle \wbar{\mathbf{g}}_t, \wbar{\boldsymbol{\lambda}}_t - \boldsymbol{\lambda}^* \rangle &\leq \frac{1}{2\alpha} \left( \|\wbar{\boldsymbol{\lambda}}_t - \boldsymbol{\lambda}^*\|^2 - \|\wbar{\boldsymbol{\lambda}}_{t+1} - \boldsymbol{\lambda}^*\|^2 \right) \nonumber \\
&\quad + \frac{\alpha}{2} K_c + C\delta_{x,t}. \nonumber
\end{align}

\paragraph*{Step 4 (Combining)} Summing the above bounds yields
\begin{align}
\Delta_t &\leq \frac{1}{2\alpha} \left( \|\wbar{\mathbf{x}}_t - \mathbf{x}^*\|^2 - \|\wbar{\mathbf{x}}_{t+1} - \mathbf{x}^*\|^2 \right) \nonumber \\
&+ \frac{1}{2\alpha} \left( \|\wbar{\boldsymbol{\lambda}}_t - \boldsymbol{\lambda}^*\|^2 - \|\wbar{\boldsymbol{\lambda}}_{t+1} - \boldsymbol{\lambda}^*\|^2 \right)  + \frac{\alpha}{2} K_c + C_E \delta_{x,t}, \nonumber
\end{align}
where $K_c$ bounds the norm of $\bar{\mathbf{p}}_t$.
\mbox{}\hfill \QED

\subsection{Proof of Theorem~\ref{thm:violation}}
The proof proceeds by relating the constraint violation to the evolution of the averaged dual variable and bounding the network disagreement.

Consider the $k$-th component of the local dual update. According to Algorithm 1, the update utilizes the consensus belief $\widehat{\mathbf{x}}_{i,t}$. With step size $\alpha = T^{-c}$:
\[
\lambda_{i,t+1}^{(k)}
=
\mathcal{P}_{[0,\Lambda_{\max}]}
\!\left(
\widehat{\lambda}_{i,t}^{(k)} + \alpha g_t^{(k)}(\widehat{\mathbf{x}}_{i,t})
\right).
\]

Averaging over all agents and using the non-negativity of the projection, we obtain:
\begin{align}
&\wbar{\boldsymbol{\lambda}}_{t+1}^{(k)}
\ge
\wbar{\boldsymbol{\lambda}}_t^{(k)}
+
\alpha \frac{1}{N}\sum_{i=1}^N g_t^{(k)}(\widehat{\mathbf{x}}_{i,t}),\nonumber\\
\Longrightarrow\quad&\frac{1}{N}\sum_{i=1}^N g_t^{(k)}(\widehat{\mathbf{x}}_{i,t})
\le
\frac{1}{\alpha}
\big(
\wbar{\boldsymbol{\lambda}}_{t+1}^{(k)} - \wbar{\boldsymbol{\lambda}}_t^{(k)}
\big).\label{eq:dual_violation_bound}
\end{align}

By the Lipschitz continuity of $g_t^{(k)}(\cdot)$,
\begin{align*}
g_t^{(k)}(\wbar{\mathbf{x}}_t)
&=
\frac{1}{N}\sum_{i=1}^N g_t^{(k)}(\widehat{\mathbf{x}}_{i,t})
+
\Big(
g_t^{(k)}(\wbar{\mathbf{x}}_t)
-
\frac{1}{N}\sum_{i=1}^N g_t^{(k)}(\widehat{\mathbf{x}}_{i,t})
\Big) \\
&\le
\frac{1}{\alpha}
\big(
\wbar{\boldsymbol{\lambda}}_{t+1}^{(k)} - \wbar{\boldsymbol{\lambda}}_t^{(k)}
\big)
+
L_g\delta_{x,t}.
\end{align*}

Summing over $t=1,\ldots,T$ and substituting $\alpha = T^{-c}$:
\[
\sum_{t=1}^T g_t^{(k)}(\widehat{\mathbf{x}}_t)
\le
\frac{1}{\alpha}\wbar{\boldsymbol{\lambda}}_{T+1}^{(k)}
+
L_g\sum_{t=1}^T \delta_{x,t}
\le
\Lambda_{\max} T^c
+
L_g\sum_{t=1}^T \delta_{x,t}.
\]

From Lemma~\ref{lemma:consensus_recur}, the consensus error satisfies the recursion:
\[
\delta_{x,t+1}
\le
\sigma \delta_{x,t} + \alpha \sqrt{N} G_d.
\]
Summing this inequality over $t=1,\ldots,T$ and leveraging the steady-state property:
\[
\sum_{t=1}^T \delta_{x,t}
\le
\frac{\sqrt{N}G_d}{1-\sigma}\, T \alpha
=
\frac{\sqrt{N}G_d}{1-\sigma}\, T^{1-c}.
\]
Combining bounds yields:
\[
\mathcal{C}_k(T) \le \Lambda_{\max} T^c + \left( \frac{L_g \sqrt{N} G_d}{1-\sigma} \right) T^{1-c}.
\]
This completes the proof. \hfill\QED

\subsection{Proof of Theorem \ref{thm:static_regret}}
The proof combines the primal--dual drift inequality in Lemma~\ref{lemma:drift}
with a bound on the network consensus error, utilizing the generic step size $\alpha = T^{-c}$.

For each $t$, recall $\wbar{\mathbf{g}}_t = \frac{1}{N}\sum_{i=1}^N \mathbf{g}_{i,t}(\mathbf{x}_{i,t})$, decompose the instantaneous regret as 
\begin{align}\label{eq:sr_decomp}
    f_t(\what{\mathbf{x}}_t) - f_t(\mathbf{x}^*) 
    &= f_t(\what{\mathbf{x}}_t) - f_t(\mathbf{x}^*) + \langle \wbar{\boldsymbol{\lambda}}_t, \wbar{\mathbf{g}}_t \rangle - \langle \wbar{\boldsymbol{\lambda}}_t, \wbar{\mathbf{g}}_t\rangle 
\end{align}
Applying Lemma~\ref{lemma:drift} with step size $\alpha=T^{-c}$ and
comparator $(\mathbf{x}^*,\mathbf{0})$, we obtain
\begin{align*}
&f_t(\widehat{\mathbf{x}}_t)-f_t(\mathbf{x}^*)
+ \langle \wbar{\boldsymbol{\lambda}}_t, \wbar{\mathbf{g}}_t \rangle \\
\le\;&
\frac{1}{2\alpha}
\Big(
\|\wbar{\mathbf{x}}_t-\mathbf{x}^*\|^2
-
\|\wbar{\mathbf{x}}_{t+1}-\mathbf{x}^*\|^2
\Big)
+
\frac{1}{2\alpha}
\Big(
\|\wbar{\boldsymbol{\lambda}}_t\|^2
-
\|\wbar{\boldsymbol{\lambda}}_{t+1}\|^2
\Big) \\
&+
\frac{\alpha}{2}(K_c+G_g^2)
+
C_E\,\delta_{x,t}.
\end{align*}
Moreover, by non-expansiveness of the dual projection,
the negative inner-product term in \eqref{eq:sr_decomp} is absorbed into
the same dual drift bound.
Summing the above inequality from $t=1$ to $T$ and using
$\|\wbar{\mathbf{x}}_t-\mathbf{x}^*\|\le D_x$ and
$\|\wbar{\boldsymbol{\lambda}}_t\|\le \Lambda_{\max}$ yield
\[
\sum_{t=1}^T
\big(f_t(\widehat{\mathbf{x}}_t)-f_t(\mathbf{x}^*)\big)
\le
\frac{D^2}{2\alpha}
+
\frac{\alpha T}{2}(K_c+G_g^2)
+
C_E \sum_{t=1}^T \delta_{x,t}.
\]
 substituting $\alpha=T^{-c}$, the first term scales as $T^c$ and the second term scales as $T^{1-c}$.

\paragraph*{Step 4: Bounding the network error}
From Lemma~\ref{lemma:consensus_recur}, the consensus error satisfies
$
\delta_{x,t+1}
\le
\sigma \delta_{x,t} + \alpha \sqrt{N} G_d .
$
Iterating this recursion and summing over $t$ yields the accumulated network error:
\[
\sum_{t=1}^T \delta_{x,t}
\le
\frac{\sqrt{N}G_d}{1-\sigma}\,T\alpha
=
\frac{\sqrt{N}G_d}{1-\sigma} T^{1-c}.
\]
Combining Step 3 and Step 4 explicitly recovers the bound in Theorem \ref{thm:static_regret}:
\[
\mathcal{R}_S(T) \le \frac{D^2}{2} T^c + \left( \frac{K_c+G_g^2}{2} + \frac{C_{\mathrm{net}}\sqrt{N}}{1-\sigma} \right) T^{1-c},
\]
where $C_{\mathrm{net}}$ combines constants $C_E$ and $G_d$. \hfill\QED

\subsection{Derivation of Inequality \eqref{eq:coupling_bound}}
In this part, we explicitly derive the upper bound on the cumulative primal consensus error under the conventional decision-sharing framework (e.g., \cite{yuan20ma-coco}), clarifying the coupling between disagreement and constraint violations.

Consider a standard distributed algorithm where agents update their decision variables $x_{i,t}$ using a gradient-based step. The accumulation of consensus error is fundamentally driven by the magnitude of these local updates.

Using standard spectral graph theory analysis, the total disagreement is bounded by the accumulation of the norms of the local update vectors:
\begin{equation} \label{eq:consensus_general}
\sum_{t=1}^T \sum_{i=1}^N \|x_{i,t} - \bar{x}_t\| \le C_{\text{net}} \sum_{t=1}^T \alpha_t \sum_{i=1}^N \|d_{i,t}\|,
\end{equation}
where $C_{\text{net}}$ is a constant related to the spectral gap of the communication network, and $d_{i,t}$ denotes the local update direction at time $t$. Since
\[
d_{i,t} = \nabla f_t(x_{i,t}) + \sum_{k=1}^m \lambda_{i,t}^{(k)} \nabla g_t^{(k)}(x_{i,t}),
\]
using the triangle inequality and the boundedness assumptions on gradients ($\|\nabla f\| \le G_f$ and $\|\nabla g\| \le G_g$), the norm is bounded by:
\begin{equation} \label{eq:norm_bound}
\|d_{i,t}\| \le G_f + G_g \|\mathbf{\lambda}_{i,t}\|.
\end{equation}
According to \cite{yuan20ma-coco}, the dual variables $\lambda_{i,t}^{(k)}$ are explicitly defined proportional to the instantaneous constraint violations, substituting this specific dual structure into \eqref{eq:norm_bound}, we obtain:
\[
\|d_{i,t}\| \le G_f + G^{'}_g \sum_{k=1}^m [g_t^{(k)}(x_{i,t})]_+.
\]
Substituting the bound for $\|d_{i,t}\|$ back into the consensus error inequality \eqref{eq:consensus_general}, and letting constants $C_1 = C_{\text{net}} G_f$ and $C_2 = C_{\text{net}} G^{'}_g$, we arrive at the coupling inequality:
\[
\sum_{t=1}^T \sum_{i=1}^N \|x_{i,t} - \bar{x}_t\| \le C_1 \sum_{t=1}^T \beta_t + C_2 \sum_{t=1}^T \beta_t \sum_{i=1}^N [g_{i,t}(x_{i,t})]_+.
\]
This confirms that in such frameworks, large constraint violations directly inflate the network disagreement, creating the coupling challenge addressed in our proposition.

\end{document}